\documentclass[12pt,reqno]{article}
mm \textheight=8.9in

\usepackage{amssymb}
\usepackage{amsmath}
\usepackage{amsthm}

\usepackage{color}
\usepackage{pb-diagram}
\newcommand{\tw}[3]{{$#1$}${\,\scriptscriptstyle {#2}}\atop\raise9pt\hbox{$\scriptstyle\tp$} ${$#3$}}
\newcommand{\st}[1]{\mbox{${\,\scriptscriptstyle {#1}}\atop\raise5.5pt\hbox{$*$}$}}
\newcommand{\btr}{\raise1.2pt\hbox{$\scriptstyle\blacktriangleright$}\hspace{2pt}}

\newcommand{\id}{\mathrm{id}}
\newcommand{\im}{\mathrm{im}\:}

\newcommand{\Tc}{\mathcal{T}}

\newcommand{\A}{\mathcal{A}}

\newcommand{\Ic}{\mathcal{I}}

\newcommand{\D}{\mathfrak{D}}

\newcommand{\Sg}{\mathfrak{S}}
\newcommand{\Jg}{\mathfrak{J}}
\newcommand{\Tg}{\mathfrak{T}}

\newcommand{\Ha}{\mathcal{H}}
\newcommand{\Ru}{\mathcal{R}}

\newcommand{\E}{\mathcal{E}}
\newcommand{\J}{\mathcal{J}}

\newcommand{\Kc}{\mathcal{K}}

\newcommand{\Nc}{\mathcal{N}}

\newcommand{\C}{\mathbb{C}}
\newcommand{\Z}{\mathbb{Z}}

\newcommand{\N}{\mathbb{N}}

\newcommand{\tp}{\otimes}

\newcommand{\U}{\mathcal{U}}
\newcommand{\F}{\mathcal{F}}
\newcommand{\ve}{\varepsilon}
\newcommand{\gm}{\gamma}

\newcommand{\tr}{\triangleright}
\newcommand{\tl}{\triangleleft}

\newcommand{\End}{\mathrm{End}}

\newcommand{\btl}{\mbox{\raise1.1pt\hbox{$\scriptstyle\blacktriangleleft$}}}
\newcommand{\ad}{\mathrm{ad}}

\newcommand{\g}{\mathfrak{g}}

\newcommand{\h}{\mathfrak{h}}

\newcommand{\si}{\sigma}
\newcommand{\al}{\alpha}

\newcommand{\be}{\begin{eqnarray}}
\newcommand{\ee}{\end{eqnarray}}

\newtheorem{thm}{Theorem}[section]
\newtheorem{propn}[thm]{Proposition}
\newtheorem{lemma}[thm]{Lemma}

\theoremstyle{definition}
\newtheorem{remark}[thm]{Remark}

\newcommand{\select}[1]{\textcolor{red}{\bf\em #1}}

\begin{document}
\title{On quantization of Semenov-Tian-Shansky Poisson bracket on
simple algebraic groups\footnote{
This research is partially supported
by
the Emmy Noether Research Institute for Mathematics,
the Minerva Foundation of Germany,  the Excellency Center "Group
Theoretic Methods in the study of Algebraic Varieties"  of the Israel
Science foundation, and by the RFBR grant no. 03-01-00593.
}
}
\author{A. Mudrov
\\
\small
\sl Dedicated to the memory of Joseph Donin
\\
}
\date{}
\maketitle
\begin{center}
{Emmy Noether Mathematics Institute, 52900 Ramat Gan,
Israel,\\
Max-Planck Institut f$\ddot{\rm u}$r Mathematik, Vivatsgasse 7, D-53111 Bonn, Germany.\\
e-mail: mudrova@macs.biu.ac.il, mudrov@mpim-bonn.mpg.de}
\end{center}
\begin{abstract}
Let $G$ be a simple complex factorizable Poisson Lie algebraic group.
Let $\U_\hbar(\g)$ be the corresponding quantum group.
We study $\U_\hbar(\g)$-equivariant quantization $\C_\hbar[G]$
of the affine coordinate ring $\C[G]$ along the Semenov-Tian-Shansky bracket.
For a simply connected group $G$
we prove
 an analog of the Kostant-Richardson theorem
stating  that $\C_\hbar[G]$ is a free module over its center.
\end{abstract}
\maketitle

\section{Introduction}
Let $G$ be a simple complex algebraic group.
Suppose $G$ is a Poisson Lie group relative to a
quasitriangular Lie bialgebra structure on $\g =\mathrm{Lie}\; G$.
Consider $G$ as a $G$-manifold with respect to the conjugation action.
In the present paper we study quantization of
a special Poisson structure on $G$
 making it a Poisson Lie $G$-manifold.
This (STS) Poisson structure is due to Semenov-Tian-Shansky. In fact, the
STS bracket  makes $G$ a Poisson Lie manifold over $\D G$, where
$\D G=G\times G$ is the Poisson Lie group corresponding to the double Lie bialgebra $\D\g\simeq \g\oplus \g$.

The affine coordinate ring $\C[G]$ can be quantized along the STS Poisson
bracket to a $\U_\hbar(\D\g)$-algebra $\C_\hbar[G]$. For $G$ connected, this
quantization can be realized as a subalgebra in $\U_\hbar(\g)$.
The algebra $\C_\hbar[G]$ is also realized as a quotient of the so called
reflection equation (RE) algebra associated with $\U_\hbar(\g)$.
For $G$ being a classical matrix group with the standard (zero weight)
DS bracket, the corresponding ideal in the RE
algebra is explicitly described.

Our main result is a quantum analog of the Kostant-Richardson theorem.
In \cite{K} Kostant proved that the algebras $\C[\g]$ and $\U(\g)$
are free modules over their subalgebras of $\g$-invariants.
Richardson generalized the case of $\C[\g]$ to the affine coordinate
ring of a semisimple complex algebraic group, \cite{R}.
Namely, if the subalgebra of invariants $I(G)$ (class functions) is polynomial,
then $\C[G]$ is a free $I(G)$-module generated by a $G$-submodule
in $\C[G]$ with finite dimensional isotypical components. We prove the analogous
statement for $\C_\hbar[G]$.

The main result of the present paper can be formulated
as follows.
\\{\bf Theorem.}
{\em
Let $G$ be a simple complex algebraic group and let $\C_\hbar[G]$ be the
$\U_\hbar(\D\g)$-equiva-riant quantization of $\C[G]$ along
the STS bracket. Then
\\
i) the subalgebra $I_\hbar(G)$ of $\U_\hbar(\g)$-invariants coincides with the center of $\C_\hbar[ G]$,
\\
ii) $I_\hbar(G)\simeq I(G)\tp \C[[\hbar]]$ as a $\C$-algebra.\\
Suppose that $I(G)$ is a polynomial algebra. Then
\\
iii) $\C_\hbar[G]$ is a free $I_\hbar(G)$-module
generated by a $\U_\hbar(\g)$-submodule $\E\subset \C_\hbar[G]$. Each isotypic component in
$\E$ is $\C[[\hbar]]$-finite.
}

Remark that for connected simply connected $G$ the algebra of invariants is a polynomial algebra generated by
the characters of fundamental representations, \cite{St}. That also is true for some non-simply connected groups,
for example, for $SO(2n+1)$.

\vspace{0.5cm}
\noindent
{\bf \large Acknowledgements.}
The author is grateful to the Max-Planck Institute for Mathematics in Bonn for hospitality and
the excellent research conditions. He thanks D. Panyushev and P. Pyatov for valuable
remarks.
\section{Quantized universal enveloping algebras}
Throughout the paper $\g$ is a simple complex Lie algebra
equipped with a quasitriangular Lie bialgebra structure. That is, we fix
a classical solution $r\in \g\tp \g$ to the Yang-Baxter equation
\be
\label{cYBE}
[r_{12},r_{13}]+[r_{12},r_{23}]+[r_{13},r_{23}]=0
\ee
and normalize it so that the symmetric part $\Omega:=\frac{1}{2}(r_{12}+r_{21})$ of $r$ is the
inverse (canonical element) of the Killing form on $\g$.
Recall that quasitriangular solutions to the equation (\ref{cYBE}) are
parameterized by combinatorial objects called Belavin-Drinfeld triples, \cite{BD}.

By $\U_\hbar(\g)$ we denote the quantization of the Lie bialgebra $(\g,r)$, \cite{Dr1,EK}.
It is a quasitriangular topological Hopf $\C[[\hbar]]$-algebra
isomorphic (as an algebra) to the space of $\U(\g)[[\hbar]]$ of formal power series in $\hbar$
with coefficients in $\U(\g)$ completed in
the $\hbar$-adic topology.

Let $\Ru\in \U_\hbar^{\hat\tp 2}(\g)$ be the quasitriangular structure (universal R-matrix) on $\U_\hbar(\g)$,
the quantization of $r\in \g^{\tp 2}$.
Consider the twisted tensor square \tw{\U_\hbar(\g)}{\Ru}{\U_\hbar(\g)} of $\U_\hbar(\g)$ constructed as follows, \cite{RS}.
The Hopf algebra
\tw{\U_\hbar(\g)}{\Ru}{\U_\hbar(\g)} is obtained by the twist of
the ordinary tensor square $\U_\hbar^{\hat\tp 2}(\g)$ by the cocycle
$\Ru_{23}\in \U_\hbar^{\hat\tp 4}(\g)$. The symbol $\hat\tp$ means completed tensor product
(in the $\hbar$-adic topology).
The diagonal embedding $\Delta \colon \U_\hbar(\g)\to \mbox{\tw{\U_\hbar(\g)}{\Ru}{\U_\hbar(\g)}}$
via comultiplication is a homomorphism of Hopf algebras. The algebra \tw{\U_\hbar(\g)}{\Ru}{\U_\hbar(\g)}
is a quantization of the double $\D\g$, which
in the simple quasitriangular (factorizable) case is isomorphic to $\g\oplus \g$
as a Lie algebra. We will use notation $\U_\hbar(\D\g)$ for \tw{\U_\hbar(\g)}{\Ru}{\U_\hbar(\g)}.

\section{Simple groups as Poisson Lie manifolds}
Given an element $\xi\in \g$ let $\xi^l$ and $\xi^r$ denote, respectively, the left- and right invariant
vector fields on $G$. Namely, we put
\be
\label{lrf}
(\xi^l f)(g)=\frac{d}{dt}f(ge^{t\xi})|_{t=0},
\quad
(\xi^r f)(g)=\frac{d}{dt}f(e^{t\xi}g)|_{t=0}
\ee
for every smooth function $f$ on $G$.

There are two important Poisson structures on $G$.
First of them, the Drinfeld-Sklyanin (DS) Poisson bracket \cite{Dr1}, is defined by the bivector field
\be
\varpi_{DS}=r^{l,l}-r^{r,r}.
\label{DSbr}
\ee
This bracket makes $G$ a Poisson Lie group, \cite{STS}.

The Semenov-Tian-Shansky (STS) Poisson structure on the group $G$ is defined
by the bivector field
\be
\varpi_{STS}=r_-^{l,l}+r_-^{r,r}-r_-^{r,l}-r_-^{l,r}
+\Omega^{l,l}-\Omega^{r,r}+\Omega^{r,l}-\Omega^{l,r}
=
r_-^{\ad,\ad} +(\Omega^{r,l}-\Omega^{l,r}).
\label{STSbr}
\ee
Here $r_-$ is the skew symmetric part $\frac{1}{2}(r_{12}-r_{21})$ of $r$.

Consider the group $G$ as a $G$-space with respect to
conjugation.
Then the  STS bracket makes $G$ a Poisson-Lie manifold
over $G$ endowed with the Drinfeld-Sklyanin bracket, \cite{STS}.

We assume $G$ to be a linear algebraic group, e.g. a subgroup
of $GL(V)$, where $V$ be a finite dimensional $G$-module.
Then $G$ is an affine variety. Its irreducible (connected) component
is an affine variety as well, \cite{VO}. Unless otherwise explicitly stated, $G$ is
assumed to be connected.

The Lie algebra $\g$ generates the left and right invariant
vector fields on $\End(V)$ defined similarly to (\ref{lrf}). Introduce a bivector field
on $\End(V)$ by the formula (\ref{STSbr}), where the superscripts $l,r$ mark the left-and right invariant vector field on $\End(V)$.
This bivector field is Poisson on the $G\times G$-invariant variety $\End(V)^\Omega$
of matrices $A\in \End(V)$ satisfying the quadratic equation $[A\tp A,\Omega]=0$.
Restriction of this Poisson structure to $G\subset \End(V)$ coincides with (\ref{STSbr}).

In the basic representation of $SL(n)$ the variety $\End(V)^\Omega$ is the entire
matrix space.
Let $G$ be an orthogonal or symplectic group
and $V$ its  basic representation with the invariant
form $B\in V\tp V$. The variety $\End(V)^\Omega$ coincides with the set of matrices
fulfilling
\be
\label{orthog''}
BX^tB^{-1}X=f^2,
\quad XBX^tB^{-1}=f^2.
\ee
Here $f$ is a numeric parameter.
The condition $f\not=0$  specifies
a principal open set in $\End(V)$, which is a group and a trivial central extension of $G$ ($f=1$).
Such an extension can be defined for an arbitrary matrix algebraic group and it will play a role
in our consideration.

\section{Quantization of the STS bracket on the group}
By \select{quantization} of a Poisson affine variety $\C[M]$ we understand
a $\C[[\hbar]]$-free $\C[[\hbar]]$-algebra $\C_\hbar[M]$ such that
$\C_\hbar[M]/\hbar \C_\hbar[M]\simeq \C[M]$. The quantization is called \select{equivariant}
if equipped with an action of a quantum group $\U_\hbar(\g)$ that is compatible with
the multiplication, namely
$$
x\tr (ab) =(x^{(1)}\tr a)(x^{(2)}\tr a)
\quad
\mbox{for all}
\quad
x\in \U_\hbar(\g)
\quad
\mbox{for all}
\quad a,b \in \C_\hbar[M]
$$
For an equivariant quantization to exist, $M$ must  be a Poisson Lie manifold
over the Poisson Lie group $G$ corresponding to the Lie bialgebra $\g$.
\subsection{Some commutative algebra}
\label{ssecQn}
In the present subsection we collect, for reader's convenience,
some standard facts about $\C[[\hbar]]$-modules that we use
in what follows.

\begin{lemma}
\label{free}
Let $E$ be a free finite $\C[[h]]$-module.
Then every  $\C[[h]]$-submodule of $E$ is finite and free.
\end{lemma}
\noindent
This assertion holds true for modules over principal ideal domains, see e.g.
\cite{Jac}.

Given an $\C[[h]]$-module $E$ we denote by $E_0$ its "classical limit",
the complex vector space $E/\hbar E$. A $\C[[\hbar]]$-linear map $\Psi\colon E\to F$
induces a $\C$-linear map $E_0\to F_0$, which will be denoted by $\Psi_0$.
\begin{lemma}
\label{Nak}
Let $E$ be a finite and $W$ an arbitrary $\C[[h]]$-modules.
A $\C[[h]]$-linear map $W\to E$ is an epimorphism if the induced map $W_0\to E_0$ is  an epimorphism
of vector spaces.
\end{lemma}
\noindent
This is a particular case of the Nakayama lemma for modules
over local rings, see e.g. \cite{GH}.

We say that a $\C[[\hbar]]$-module $E$ has no torsion (is torsion free) if $\hbar x=0 \Rightarrow x=0$ for $x\in E$.
\begin{lemma}
\label{torsion_free}
A finitely generated $\C[[h]]$-module is free if it is torsion free.
\end{lemma}
\noindent
The latter assertion easily follows from the Nakayama lemma.

\begin{lemma}
\label{sub-fac}
Every submodule and quotient module of a finite $\C[[\hbar]]$-module is finite.
\end{lemma}
\noindent
This statement is obvious for quotient modules. For submodules, it follows from
Lemma \ref{free}.

\begin{lemma}
\label{auto}
Let $\Psi\colon E\to F$ be a morphism of free finite $\C[[\hbar]]$-modules
such that the induced map $\Psi_0\colon E_0\to F_0$ is an isomorphism of $\C$-vector spaces.
Then $\Psi$ is an isomorphism.
\end{lemma}
\noindent
Using Lemma \ref{free}, the latter assertion can be reduced to the case $E=F$
and $\Psi$ being an endomorphism of $E$. An endomorphism of a free
module is invertible if
and only if its residue $\mod \hbar$ is invertible.

\begin{lemma}
\label{injective}
Let $\Psi\colon E\to F$ be a morphism of a $\C[[\hbar]]$-modules.
Suppose that $E$ is finite, $F$ is torsion free, and $\Psi_0\colon E_0\to F_0$
is injective. Then $E$ is free, and $\Psi$ is injective.
\end{lemma}
\begin{proof}
First let us prove that $\Psi$ is embedding assuming $E$ to be free.
In this case the image $\im \Psi$ is finite and has no torsion.
Therefore it is free, by Lemma \ref{torsion_free}.
The map $\Psi_0$ factors through the composition $E_0\to (\im \Psi)_0\to F_0$,
and the left arrow is surjective by construction.
Since $\Psi_0$ is injective, the map $E_0\to (\im \Psi)_0$ is also injective and hence
an isomorphism, by Lemma \ref{auto}. Therefore $E\simeq \im \Psi$.

Now let $E$ be arbitrary and let $\{e_i\}$ be a
set of generators such that their projections $\mod \hbar$ form a base in $E_0$.
Such generators do exist in view of the Nakayama lemma. Let $\hat E$
be the $\C[[\hbar]]$-free covering of $E$ generated by $\{e_i\}$. The composite map
$\hat \Psi \colon \hat E\to E\to F$ satisfies the hypothesis of the lemma
with free $\hat E$. We conclude that $\hat\Psi$ is injective.
This implies that $E=\hat E$, i. e. $E$ is free, and that $\Psi$ is injective.
\end{proof}
\subsection{Quantization of the DS and STS brackets}
\label{ssecQDSSTS}
Let $G$ be a simple complex algebraic group and $(V,\rho)$ its faithful representation.
The affine ring $\C[G]$ is realized as a quotient of $\C[\End(V)]$ by
an ideal generated by a finite system of polynomials $\{p_i\}$.
We do not require $G$ to be connected but assume
$G\subset \End(V)^\Omega$. Recall that  $\Omega\in \g\tp \g$ stands for the split Casimir.

Let $G^\sharp$ denote the smooth affine variety $G\times \C^*$, where $\C^*$ is the multiplicative
group of the field $\C$. The variety $G^\sharp$ is an algebraic group, however we will not use this
fact until Section \ref{gen}.

The affine coordinate ring
$\C[G^\sharp]$ is isomorphic to the tensor product $\C[G]\tp \C[f,f^{-1}]$.
It can be realized as the quotient of $\C[\End(V)]\tp \C[f,f^{-1}]$
by the ideal $(p_i^\sharp)$, where $p_i^\sharp(f,X)=f^{k_i}p_i(f^{-1}X)$ and $k_i$ is the
degree of the polynomial $p_i$.

The algebra $\C[\End(V)]\tp \C[f,f^{-1}]$ is equipped with a $\Z$-grading by setting
$\deg \End^*(V)=1$, $\deg f=1$, and $\deg f^{-1}=-1$. The polynomials $\{p_i^\sharp\}$
are homogeneous, hence $\C[G^\sharp]$ is a $\Z$-graded algebra.
Let us select in $\C[G^\sharp]$  the subalgebra
which is the quotient of $\C[\End(V)][f]$ by the ideal $(p_i^\sharp)$.
This subalgebra is identified with the affine ring of the Zariski
closure $\bar G^\sharp$
in $\End(V)\times \C$. It is graded, with finite dimensional
homogeneous components.
Clearly $\C[G^\sharp]$ is generated by $\C[\bar G^\sharp]$ over $\C[f^{-1}]$.

Define a two-sided $G$-action on $\C[G^\sharp]$ by setting it trivial on $\C[f,f^{-1}]$.
This makes $\C[G^\sharp]$ a $\U(\g)$-bimodule algebra. The  action preserves grading
and preserves the subalgebra  $\C[\bar G^\sharp]$.
The DS and STS brackets (\ref{DSbr}) and (\ref{STSbr}) are naturally defined on $\C[G^\sharp]$ and  $\C[\bar G^\sharp]$ via the right and left $\g$-actions
on $\C[G^\sharp]$ and $\C[\bar G^\sharp]$ . They make both
$\C[G^\sharp]$ and $\C[\bar G^\sharp]$ Poisson Lie algebras over
the Lie bialgebras $\g_{op}\oplus \g$ and $\D\g$, correspondingly.
The Poisson Lie manifolds $G_{DS}$ and $G_{STS}$ are sub-manifolds in $G_{DS}^\sharp$
and $G_{STS}^\sharp$ (as well as in $\bar G_{DS}^\sharp$ and $\bar G_{STS}^\sharp$)
defined by the equation $f=1$.

Recall the Takhtajan quantization of the DS Poisson structure on $G$, \cite{T}.
Consider the quasitriangular quasi-Hopf algebra $\bigl(\U(\g)[[\hbar]],\Phi,\Ru_0\bigr)$,
where $\U(\g)[[\hbar]]$ is equipped with the standard comultiplication,
$\Phi$ is a $\g$-invariant  associator, and $\Ru_0=e^{\frac{\hbar}{2}\Omega}$ is the universal R-matrix.
Since $\Phi$ and $\Ru_0$ are $G$-invariant, $\C[G]\tp \C[[\hbar]]$ is a commutative algebra in the quasi-tensor
category of $\U(\g)_{op}[[\hbar]]\hat \tp\U(\g)[[\hbar]]$- modules.
 The latter
is a quasi-Hopf algebra with the associator $(\Phi^{-1})'\Phi''$ and the universal R-matrix
$(\Ru^{-1}_0)'\Ru_0''$, \cite{Dr3}. Here the prime is relative to the $\U(\g)_{op}[[\hbar]]$-factor
while the double prime to the $\U(\g)[[\hbar]]$-factor.

Let $\J\in \U(\g)^{\hat \tp 2}[[\hbar]]$ be a twist making $\U(\g)[[\hbar]]$
the quasitriangular Hopf algebra $\U_\hbar(\g)$. Then $(\J^{-1})'\J''$ converts
$\U(\g)_{op}[[\hbar]]\hat\tp \U(\g)[[\hbar]]$ into the Hopf algebra
$\U_\hbar(\g)_{op}\hat\tp\U_\hbar(\g)$. Applied to $\C[G]\tp \C[[\hbar]]$, this twist
makes it a $\U_\hbar(\g)_{op}\hat\tp \U_\hbar(\g)$-module algebra, $\C_\hbar[G_{DS}]$.
This algebra is commutative in the category of $\U_\hbar(\g)$-bimodules. It is a quantization
of the DS-Poisson Lie bracket on $G$.

The above quantization extends to the algebras $\C_\hbar[G^\sharp_{DS}]$ and $\C_\hbar[\bar G^\sharp_{DS}]$;
the construction is literally the same. Since the two-sided
action of $\g$ preserves the grading, the algebras $\C_\hbar[G^\sharp_{DS}]$ and $\C_\hbar[\bar G^\sharp_{DS}]$
are $\Z$-graded. The algebra $\C_\hbar[G_{DS}]$ is obtained
from $\C_\hbar[G^\sharp_{DS}]$ or from $\C_\hbar[\bar G^\sharp_{DS}]$  as the quotient by the ideal $(f-1)$.

Now consider  $\C_\hbar[G_{DS}]$, $\C_\hbar[G^\sharp_{DS}]$, and $\C_\hbar[\bar G^\sharp_{DS}]$
as $\U_\hbar(\g)^{op}\hat\tp\U_\hbar(\g)$-algebras,
using identification between $\U_\hbar(\g)^{op}$ and $\U_\hbar(\g)_{op}$ via the antipode.
Perform the twist  from $\U_\hbar(\g)^{op}\hat\tp\U_\hbar(\g)$ to
$\U_\hbar(\D\g)$
and the corresponding transformation of the algebras
$\C_\hbar[G_{DS}]$,  $\C_\hbar[G^\sharp_{DS}]$, and $\C_\hbar[\bar G^\sharp_{DS}]$. The resulting algebras
$\C_\hbar[G_{STS}]$, $\C_\hbar[G^\sharp_{STS}]$, and $\C_\hbar[\bar G^\sharp_{STS}]$ are $\U_\hbar(\D\g)$-equivariant quantizations
along the STS bracket, \cite{DM}. They are commutative in the braided category of $\U_\hbar(\D\g)$-modules.

The algebras $\C_\hbar[G^\sharp_{STS}]$
and $\C_\hbar[\bar G^\sharp_{STS}]$ are $\Z$-graded
and $\C_\hbar[G^\sharp_{STS}]=\C_\hbar[\bar G_{STS}][f^{-1}]$.
The homogeneous components in $\C_\hbar[\bar G_{STS}^\sharp]$ are
$\C[[\hbar]]$-finite and vanish for negative degrees.
The algebra $\C_\hbar[G_{STS}]$ is obtained from $\C_\hbar[G^\sharp_{STS}]$ (or from $\C_\hbar[\bar G^\sharp_{STS}]$) by factoring out
the ideal $(f-1)$.

\section{The algebra $\C_\hbar[G_{STS}]$ as a module over its center}
In the present section $G$ is connected and $\C_\hbar[G]$ stands for $\C_\hbar[G_{STS}]$, that is, for the
$\U_\hbar(\D\g)$-equivariant quantization of $\C[G]$ along the STS bracket.
The action of $\U_\hbar(\g)$ is induced by the diagonal embedding
$\Delta\colon \U_\hbar(\g)\to \U_\hbar(\D\g)$ and can be expressed through
the left and right coregular actions of $\U_\hbar(\g)$ on $\C_\hbar[G_{DS}]$ as
$$
x (a)= x^{(2)}\tr a\tl \gm( x^{(1)}).
$$
Here $\gm$ stands for the antipode in $\U_\hbar(\g)$ and
the actions are defined by $\xi\tr a= \xi^l(a)$, and  $a\tl \xi= \xi^r(a)$
for $\xi\in \g$, cf. (\ref{lrf}).
We use that fact that $\C_\hbar[G_{DS}]$ and $\C_\hbar[G_{STS}]$ coincide as
$\U_\hbar(\g)$-bimodules (but not algebras) and the $\U_\hbar(\g)$-actions
is the actions of $\U(\g)[[\hbar]]$.
\begin{propn}
\label{center}
Let $G$ be a simple complex algebraic group equipped with the STS bracket. Let $\g$ be
its Lie bialgebra, $\D\g$ the double of $\g$, and let $\C_\hbar[G]$ be
the $\U_\hbar(\D\g)$-equivariant quantization of the affine ring $\C[G]$ along the STS bracket.
Then the subalgebra $I_\hbar(G)$ of $\U_\hbar(\g)$-invariants in $\C_\hbar[G]$ coincides with the
center.
\end{propn}
\begin{proof}
The statement holds true for $\bar G^\sharp$ too. Let us prove it for $\bar G^\sharp$ first.
The case of $G$ will be obtained
by factoring out the ideal $(f-1)$.

The subalgebra $I_\hbar(\bar G^\sharp)$ lies in the center of $\C_\hbar[\bar G^\sharp]$.
Indeed, let $\hat \Ru$ be the universal R-matrix of $\U_\hbar(\D\g)$.
It is expressed through the universal R-matrix $\Ru \in \U_\hbar^{\hat \tp 2}(\g)$
by $\hat\Ru=\Ru_{41}^{-1}\Ru_{31}^{-1}\Ru_{24}\Ru_{23}$, therefore
$\hat \Ru\in \U_\hbar(\D\g)\hat \tp \U_\hbar(\g)$.
The algebra $\C_\hbar[\bar G^\sharp]$ is commutative in the category of $\U_\hbar(\D\g)$-modules,
hence
$(\hat \Ru_2\tr a)(\hat \Ru_1\tr b)=ba$
for any  $a,b\in \C_\hbar[\bar G^\sharp]$.
Hence $ab=ba$ for $a\in I_\hbar(\bar G^\sharp)$.

Conversely, suppose that $ab =ba$ for some $a$ and all $b\in \C_\hbar[\bar G^\sharp]$.
Present $a$ as $a=a_0+O(\hbar)$,
where $a_0\in \C[\bar G^\sharp]$.
We have $0=\hbar \varpi_{STS}(a_0,b)+O(\hbar^2)$ and therefore $\varpi(a_0,b)=0$.
The Poisson bivector field $\varpi_{STS}$ is induced by the classical r-matrix of
the double $\xi^i\tp \xi_i\in (\D\g)^{\tp 2}$. The element $\xi\in \g^*$ acts on $\bar G^\sharp$ by vector field
$r_-(\xi)^l-r_-(\xi)^r+\frac{1}{2}\bigl(\Omega(\xi)^l+\Omega(\xi)^r\bigr)$ (here we consider the
elements of $\g\tp \g$ as operators $\g^*\to \g$ by paring with the first tensor component).
Let $e$ be the identity of the group $G$.
At every point $(e\tp c)\in G\times \C^*=G^\sharp\subset \bar G^\sharp$ this vector field equals $\Omega(\xi)$.
Since the Killing form is non-degenerate, $\zeta\tr a_0=0$ for all $\zeta\in \g$
in an open set  in $\bar G^\sharp$ (in Euclidean topology). Therefore $\zeta\tr a_0=0$ for all $\zeta\in \g$ and $a_0$ is $\g$-invariant.

We can assume that $a$ is homogeneous with respect to the grading
in $I_\hbar(\bar G^\sharp)$.
Let $a'_0$ be $\U_\hbar(\g)$-invariant element such that $a'_0=a_0\mod \hbar$.
We can choose $a'_0$ of the same degree as $a$ (in fact, we can take $a'_0=a_0\tl \theta^{-\frac{1}{2}}$,
see the proof of Proposition \ref{non-def}).
Then $a-a'_0$ is central and divided by $\hbar$.
Acting by induction, we present $a$ as a sum $a=\sum_{\ell=0}^\infty \hbar^\ell a'_\ell$,
where each summand is $\U_\hbar(\g)$-invariant. Since all $a'_\ell$ have the same degree,
they lie in a finite $\C[[\hbar]]$-module. Hence the above sum converges in the $\hbar$-adic topology.
\end{proof}
An immediate corollary of Proposition \ref{center} is the analogous statement for $G^\sharp$ and $\bar G^\sharp$.

The following proposition asserts that the subalgebra of invariants in $\C_\hbar[G]$ is not quantized.
\begin{propn}
\label{non-def}
Let $\C_\hbar[G]$ be the $\U_\hbar(\D\g)$-equivariant quantization of the STS bracket on $G$.
Then $I_\hbar(G)$ is isomorphic to $I(G)\tp \C[[\hbar]]$ as a $\C$-algebra.
\end{propn}
\begin{proof}
Consider two subspaces $\Ic_1$ and $\Ic_2$ in $\A=\C_\hbar[G_{DS}]$ defined by the following
conditions:
\be
\label{2inv}
\Ic_1=\{a\in \A\colon x\tr a = a\tl x, \forall x\in \U_\hbar(\g)\},
\quad
\Ic_2=\{a\in \A\colon x\tr a = a\tl \gm^2(x), \forall x\in \U_\hbar(\g)\}.
\ee
Since $\A$ is an $\U_\hbar(\g)$-bimodule algebra, and
the square antipode $\gm^2$ is a Hopf algebra automorphism of $\U_\hbar(\g)$,
both $\Ic_1$ and $\Ic_2$ are subalgebras in $\A$.
The algebra $\Ic_1$ is isomorphic to $I(G)\tp \C[[\hbar]]$ as a $\C$-algebra.
This readily follows from the Takhtajan construction of $\C_\hbar[G_{DS}]$ rendered
in Subsection \ref{ssecQDSSTS}.

Let us show that the algebra $\Ic_2$ is isomorphic to $\Ic_1$.
Indeed, the fourth power of the antipode in $\U_\hbar(\g)$
is implemented by the  similarity transformation with a group-like element
$\theta\in \U_\hbar(\g)$, \cite{Dr2}. This element has a group-like square root
$\theta^{\frac{1}{2}}=e^{\frac{1}{2}\ln \theta}\in \U_\hbar(\g)$. The logarithm is well defined,
because $\theta=1+O(\hbar)$.
In the case of the Drinfeld-Jimbo or standard quantization of $\U(\g)$
the element $\theta^{\frac{1}{2}}$ belongs to $\U_\hbar(\h)$, where $\h\subset \g$ is the Cartan subalgebra.
The map $a\mapsto a\tl \theta^{-\frac{1}{2}}$ is an automorphism of $\A$, and this
automorphism sends $\Ic_1$ to $\Ic_2$.

Thus we have proven that $\Ic_2$ is isomorphic to $I(G)\tp \C[[\hbar]]$ as a $\C$-algebra.
Consider the RE twist converting $\C_\hbar[G_{DS}]$ into $\C_\hbar[G_{STS}]$.
This twist relates multiplications by the formula (\ref{tw_pr}), where
 $\Tc$ should be replaced by $\C_\hbar[G_{DS}]$ and $\Kc$ by $\C_\hbar[G_{STS}]$.
It is straightforward to see that these multiplications coincide on $\Ic_2$.
\end{proof}
\begin{remark}
\label{left}
In the proof of Proposition \ref{non-def}, we used the observation that the
 multiplications in $\C_\hbar[G_{DS}]$ and $\C_\hbar[G_{STS}]$
coincide on $\Ic_2$. In fact, formula (\ref{tw_pr}) implies that  $\C_\hbar[G_{DS}]$ and $\C_\hbar[G_{STS}]$
are the same as left $\Ic_2$-modules. Therefore the structure of left $\Ic_1$-module on
$\C_\hbar[G_{DS}]$ is the same as the structure of $\Ic_\hbar(G)$-module on $\C_\hbar[G_{STS}]$.
This assertion also holds for $G^\sharp$ and $\bar G^\sharp$.
\end{remark}

Let $T$ be the maximal torus in $G$.
Then $I(G)\simeq \C[T]^W$, where $W$ is the Weyl group, \cite{St}.
Suppose that the subalgebra of invariants in $\C[G]$ is polynomial. For example, that is the case when $G$ is simply
connected; then $\C[T]^W$ is generated by characters of the fundamental representations \cite{St}.
Under the above assumption, the algebra $\C[G]$ is a free module over $I(G)$, \cite{R}.
There exists a $G$-submodule $\E_0\subset \C[G]$ such that the multiplication
map $I(G)\tp \E_0\to \C[G]$ gives an isomorphism of vector spaces.
Each isotypic component in $\E_0$ has finite multiplicity. We will establish
the quantum analog of this fact.

\begin{thm}
\label{G=IE}
Let $\C_\hbar[G]$ be the $\U_\hbar(\D\g)$-equivariant quantization of the STS bracket on $G$.
Suppose that the subalgebra $I(G)$ of $\g$-invariants is a polynomial algebra.
Then
\\
i) $\C_\hbar[G]$ is a free $I_\hbar(G)$-module generated by
a $\U_\hbar(\g)$-submodule $\E\subset \C_\hbar[G]$.
\\ii)
each isotypic component in $\E$ is $\C[[\hbar]]$-finite.
\end{thm}
\begin{proof}
Let $\E_0$ be the $\U(\g)$-module generating $\C[G]$ over $I(G)$.
Naturally considered as a subspace in $\C[ G^\sharp]$, it
obviously generates $\C[ G^\sharp]$ over $I(G^\sharp)$.
Using invertibility of $f$, we can make  every isotypic component of $\E_0$ homogeneous and regard $\E_0$ as a
graded submodule in $\C[\bar G^\sharp]$.

Put $\E=\E_0\tp \C[[\hbar]]$.
Let $V_0$ be a simple finite dimensional $\g$-module and $V=V_0\tp \C[[\hbar]]$
the corresponding $\U_\hbar(\g)$-module. Let
$(\E_0)_{V_0}$ denote the isotypic component of
$\E_0$. The isotypic component $\C[G^\sharp]_V$ is isomorphic
to $I(G^\sharp)\tp (\E_0)_{V_0}\tp \C[[\hbar]]$, as a $\U_\hbar(\g)$-module.

Let $\tilde m$ denote the multiplication in $\C_\hbar[G^\sharp]$.
The map
\be
\label{free_map}
\tilde m \colon I_\hbar(G^\sharp)\tp_{\C[[\hbar]]} \E_V \to \C_\hbar[G^\sharp]_V
\ee
is $\U_\hbar(\g)$-equivariant and respects grading.
Let the superscript $(k)$ denote the homogeneous component of degree $k$.
The map (\ref{free_map})
induces $\U_\hbar(\g)$-equivariant maps
\be
\label{restricted}
\oplus_{i+j=k} I_\hbar(G^\sharp)^{(i)}\tp_{\C[[\hbar]]} \E_V^{(j)} \to \C_\hbar[ \bar G^\sharp]_V^{(k)}
,
\quad
I_\hbar(\bar G^\sharp)^{(k)}\tp_{\C[[\hbar]]} \E_V \to \C_\hbar[\bar G^\sharp]_V\subset \C_\hbar[G^\sharp]_V.
\ee
The left map has a $\C[[\hbar]]$-finite target, while the right one has a $\C[[\hbar]]$-finite source.
All the $\C[[\hbar]]$-modules in (\ref{restricted}) are free.
Modulo $\hbar$, the left map is surjective, and the right one injective.
Therefore they are surjective and injective, respectively, by Lemmas \ref{Nak} and \ref{injective}.
Since $\C_\hbar[G^\sharp]=\C_\hbar[\bar G^\sharp][f^{-1}]$ and $I_\hbar[G^\sharp]=I_\hbar[\bar G^\sharp][f^{-1}]$, this
immediately implies that the map (\ref{free_map}) is surjective and injective and hence an isomorphism.

Now recall that $I_\hbar(G^\sharp)$ is  isomorphic to $I_\hbar(G)[f,f^{-1}]$.
Taking quotient by the ideal $(f-1)$ proves the theorem for $G$.
\end{proof}
\section{Quantization in terms of generators and relations}
\label{gen}
In this section we describe the quantization of $\C[G]$ along the DS and STS brackets in terms of generators and relations
for $G$ being
a classical matrix  group. We give a detailed consideration to the DS-case. The case of STS is treated
similarly, upon obvious modifications. Alternatively, the defining ideal
$\C_\hbar[G_{STS}]$ can be derived from the ideal of $\C_\hbar[G_{DS}]$ using
Proposition \ref{propn_ap} and the twist-equivalence between $\C_\hbar[G_{DS}]$ and
$\C_\hbar[G_{STS}]$.

Function algebras on quantum classical matrix groups from the classical series
were defined in terms of generators and relations  in
\cite{FRT}. Here we prove that the algebras of \cite{FRT} are included in  flat $\C[[\hbar]]$-algebras,
$\C_\hbar[G_{DS}]$.

\subsection{FRT and RE algebras}
\label{ssecqSTS}
In this subsection we recall the definition of the FRT and RE algebras, \cite{FRT,KSkl}.

Let $V_0$ be the basic representation of $G$ and let $V$ be the corresponding
$\U_\hbar(\g)$-module. Let $R$ denote the
image of the universal R-matrix of $\U_\hbar(\g)$ in $\End(V^{\tp 2})$. Put  $N:=\dim V_0$,

The FRT algebra $\Tc$ is generated by the matrix coefficients $\{T^i_j\}\subset \End^*(V)$
 subject to the relations
\be
R T_1 T_2 = T_2 T_1 R,
\label{RTT}
\ee
where $T=||T^i_j||$. So $\Tc$ is the quotient of the free algebra $\C[[\hbar]]\langle T^i_j\rangle$.
The latter is a $\U_\hbar(\g)$-bimodule algebra, the two-sided action being extended
from the two-sided action on $\End^*(V)$.
The ideal (\ref{RTT}) is invariant, so $\Tc$ is also a $\U_\hbar(\g)$-bimodule algebra.
It is $\Z$-graded with $\deg \End^*(V)=1$, and the grading is equivariant with respect to the
two-sided $\U_\hbar(\g)$-action.

The RE algebra $\Kc$ is also generated by the matrix coefficients of the  basic representation,
this time denoted by $K^i_j$.
Let $K=||K^i_j||$ be the matrix of the generators.
 The RE algebra $\Kc$ is the quotient of the free algebra $\C[[\hbar]]\langle K^i_j\rangle$ by the ideal
 generated by the relations
\be
R_{21}K_1 R_{12} K_2 = K_2 R_{21} K_1 R_{12}.
\label{RE}
\ee
The algebra $\Kc$ is a $\U_\hbar(\D\g)$-module algebra, \cite{DM}. It is
$\Z$-graded, and the grading is invariant with respect to the
$\U_\hbar(\D\g)$-action.

Recall from \cite{DM} that the RE twist of the Hopf algebra $\U_\hbar^{op}(\g)\hat\tp\U_\hbar(\g)$
to the twisted tensor square $\U_\hbar(\D\g)$,
converts the algebra $\Tc$ to $\Kc$ (cf. also Subsection \ref{ssecQDSSTS}).

\subsection{Algebra $\C_\hbar[G_{DS}]$ in generators and relations.}
From now one we assume the standard (zero weight) Lie bialgebra structure on $\g$.
In this section we describe the algebra $\C_{\hbar}[G]=\C_\hbar[G_{DS}]$ in terms of
generators and relations.

We will use the group structure on $G^\sharp$, which is the trivial central extension of $G$.
For  $G$ orthogonal and symplectic, $G^\sharp$ is defined by equation (\ref{orthog''})
with $f\not =0$.
The  basic representation of $G$ on $V_0$ naturally extends
to a representation of of $G^\sharp$ on $V_0\oplus \C$, since
the subgroup $\C^*\subset G$ acts on the module $V_0$ by the delations.
The indeterminant $f$ is the matrix coefficient of the one dimensional representation of $\C^*$.

Suppose $f\not =0$. The group $G^\sharp$ can be identified with
the $G^\sharp\times G^\sharp$-orbit in $\End(V_0\oplus \C)$,
which for $G$ orthogonal and symplectic
is defined by the equation
\be
\label{G+}
&B_0T^tB^{-1}_0T=f^2,
\quad TB_0T^tB^{-1}_0=f^2,
\\
\label{GL}
&\det(T)=f^N
\ee
and by (\ref{GL}) for $G=SL(n)$.
The element $B_0\in V_0\tp V_0$ in equation (\ref{G+}) is the classical invariant of the (orthogonal or symplectic)
 group $G$.

Clearly the ideals in $\C[\End(V_0)][f]$ generated by (\ref{G+}) and by (\ref{GL}) are radical.
This is obvious for the $G=SL(n)$ and follows from \cite{We} for
$G$ orthogonal and symplectic.
The corresponding quotients of $\C[\End(V_0)][f]$ are the affine coordinate rings of $\bar G^\sharp$.

Recall from \cite{FRT} and \cite{F} that there exists a central group-like
two-sided
$\U_\hbar\bigl(sl(n)\bigr)$-invariant $\det_q(T)\in \Tc$ of degree $n$ such
that $\det_q(T)=\det(T)$ modulo $\hbar$.
For $G$ orthogonal and symplectic let $B$ denote the $\U_\hbar(\g)$-invariant element from $V\tp V$,
see \cite{FRT}.
\begin{propn}
\label{sharp}
Let $G$ be a classical unimodular matrix group.
The $\U_\hbar(\g)_{op}\tp \U_\hbar(\g)$-equivariant quantization $\C_\hbar[\bar G^\sharp]$
can be realized as the quotient of $\Tc[f]$ by the ideal of relations
\be
\label{GL'}
{\det}_q(T)=f^N
\ee
and, for $\g$ orthogonal or symplectic,
\be
\label{orthog'}
BT^tB^{-1}T=f^2,
\quad TBT^tB^{-1}=f^2.
\ee
The quantization $\C_\hbar[G]$ is obtained from $\C_\hbar[ \bar G^\sharp]$ by
factoring out the ideal $(f-1)$.
\end{propn}
\begin{proof}
Denote by $\Sg$ the algebra $\Tc[f]$, by $\Tg$ the algebra $\C_\hbar[G^\sharp]$,
and by $\Jg$ the ideal in $\Tc[f]$ generated by the relations (\ref{orthog'})
and (\ref{GL'}), depending on the type of $G$.
The algebras $\Sg$, $\Tg$, and $\Jg$ are graded, and the grading is
$\U_\hbar(\g^\sharp)$-compatible. Note that homogeneous components in $\Sg$ and
hence in $\Jg$ are $\C[[\hbar]]$-finite. There is an obvious bialgebra structure
on $\Sg$, with $f$ being group-like.

The Takhtajan construction of the quantization, see Subsection \ref{ssecQDSSTS}, implies that
the evaluation at the identity $\ve\colon a\mapsto a(e)$ is a character of
the algebra $\Tg$. Define a pairing
between $\Tg$ and  $\U_\hbar(\g^\sharp)$ setting $\langle a,x\rangle:=\ve(x\tr a)=
\ve(a\tl x)$. This pairing is non-degenerate, because $G$ is connected.

The matrix coefficients of the  basic representation are naturally considered as the elements
of $\Tg$. They  satisfy the RTT relation, because $\Tg$ is commutative in the category of
$\U_\hbar(\g^\sharp)$-bimodules. This defines an equivariant algebra homomorphism
$\Psi\colon \Sg \to \Tg$. Clearly the composition map $\ve\circ \Psi$ coincides with
the counit of the bialgebra $\Sg$.
From this we conclude that the invariant ideal $\Jg$ is
annihilated by $\ve$ (the counit gives $1$ on group-like elements,
including $\det_q(T)$ and $f$).
Therefore $\Jg$ annihilates $\U_\hbar(\g^\sharp)$
through the pairing $\langle .,.\rangle$.
Since this pairing is non-degenerate, the ideal $\Jg$ lies in the kernel of $\Psi$.

The homomorphism $\Psi$ preserves grading and it is identical
on $\End^*(V)\oplus \C[[\hbar]] f$.
As the image of $\Psi$ is $\C[[\hbar]]$-free, we have the direct sum
decomposition $\Sg=\ker \Psi\oplus \im \Psi$ of $\C[[\hbar]]$-modules.
Therefore $(\im \Psi)_0$ is embedded in $\Sg_0$.
Let us show that $\Jg=\ker \Psi$.
Since both ideals are graded and the homogeneous components are finite, it
suffices to show that the map $\Jg_0\to (\ker \Psi)_0$ induced by the embedding
$\Jg\hookrightarrow \ker \Psi$ is surjective.
Then we can apply the Nakayama lemma to each homogeneous component.

Denote by $\Jg_0^\flat$ the image of $\Jg_0$ in $(\ker \Psi)_0\subset \Sg_0$.
This is a $G^\sharp\times G^\sharp$-invariant ideal, and it is easy to show
that $\Jg_0^\flat$ contains no positive integer powers of $f$.
On the other hand,
the defining ideal $\Nc( \bar G^\sharp)\subset \Sg_0$ is maximal among such ideals,
and it
lies in $\Jg_0^\flat$. Therefore $\Nc( \bar G^\sharp)=\Jg_0^\flat =(\ker \Psi)_0$,
since otherwise the map $\Psi_0$ and hence $\Psi$ would be zero.
This proves $\Jg=\ker \Psi$.
Another consequence is that $\im \Psi$ is a quantization of $\C[ \bar G^\sharp]$ that lies in
$\C_\hbar[ \bar G^\sharp]$. Hence it coincides with $\C_\hbar[ \bar G^\sharp]$,
because that is so in the classical limit.

The quantization $\C_\hbar[G^\sharp]$ is isomorphic to $\C_\hbar[ \bar G^\sharp][f^{-1}]$,
as easily follows from the Takhtajan construction.
Therefore $\C_\hbar[G^\sharp]$ is realized as the quotient
of the algebra $\Tc[f,f^{-1}]$ by the ideal of the relations (\ref{GL'}), (\ref{orthog'}).
On the other hand, $\C_\hbar[G^\sharp]$ is a free module over
$\C[[\hbar]][f,f^{-1}]$. The quotient of $\C_\hbar[G^\sharp]$ by the ideal $(f-1)$
is $\C[[\hbar]]$-free and thus is a quantization of $\C[G]$.
\end{proof}

For orthogonal and symplectic $\g$ a stronger assertion can be proven.
Now let  $G$ be  $O(N)$ or $Sp(n)$.
\begin{propn}
\label{sharp_so}
The $\U_\hbar(\g)_{op}\tp \U_\hbar(\g)$-equivariant quantization $\C_\hbar[\bar G^\sharp]$
can be realized as the quotient of $\Tc[f]$ by the ideal of relations
(\ref{orthog'}). The quantization $\C_\hbar[G]$ is obtained from $\C_\hbar[ \bar G^\sharp]$ by
factoring out the ideal $(f-1)$.
\end{propn}
\begin{proof}
For the group $Sp(n)$ the statement is, in fact, already proven, because
$Sp(n)$ defined by (\ref{G+}) with $f=1$ is unimodular and the condition $\det_q=1$ is excessive.
Thus let us focus on the
orthogonal case.

The group $\Z_2$ acts on $\U_\hbar(\g)$ by Hopf algebra automorphisms.
This action is trivial for $\g=so(2n+1)$ and induced by the flip $\si$ of the simple roots $\al_{n-1}$ and $\al_n$
(the automorphism of the Dynkin diagram)
in the quantum Chevalley basis
for $\g=so(2n)$.
Consider the smash product $\Z_2\ltimes\U_\hbar(\g)$ with the natural structure of Hopf algebra,
a deformation of the Hopf algebra $\Z_2\ltimes\U(\g)$.
Let $\{e^i_j\}$ be the standard matrix base of
$\End(V)$.
The representation of $\U_\hbar(\g)$ on $V$ extends to a representation of
$\Z_2\ltimes\U_\hbar(\g)$ by assigning $\si\mapsto -1$ for  $\g=so(2n+1)$ and
$\si\mapsto  1-e^{n}_n-e_{n+1}^{n+1}+e^{n+1}_n+e_{n+1}^n$ for $\g=so(2n)$ (in the
realization of \cite{FRT}). Note that these matrices
are characters of the algebra $\Tc$ with $\det=-1$ in the classical limit).

Now repeat the proof of  Proposition \ref{sharp} replacing
 $\U_\hbar(\g^\sharp )$ by $\Z_2\ltimes\U_\hbar(\g^\sharp )$.
\end{proof}
It is crucial in the proofs of Propositions \ref{sharp_so} and \ref{sharp_so} that the defining ideal
of the group $G^\sharp$ should be maximal invariant.
Without the condition $\det=f^N$, the group $G^\sharp$ has two connected components, which
are orbits of the two-sided action $\U(\g^\sharp)$. Hence the defining ideal
of $G^\sharp$ is not maximal among the invariant proper ideals in $\C[G^\sharp]$. It becomes
so if we extend the symmetries and consider the algebra $\Z_2\ltimes \U(\g)$,
instead of $\U(\g)$.
\subsection{Algebra $\C_\hbar[G_{STS}]$ in generators and relations}

Under the twist from $\U_\hbar^{op}(\g)\hat\tp \U_\hbar(\g)$ to $\U_\hbar(\D\g)$,
the defining relations of $\C_\hbar[G_{STS}]$ in $\Tc$ transform to certain relations in the RE algebra
$\Kc$ and generate a $\U_\hbar(\D\g)$-invariant ideal
in $\Kc$, see Appendix \ref{apA}. Let us compute this ideal.

The multiplications in $\Tc$ and $\Kc$ are related by the formula (see \cite{DM})
\be
\label{tw_pr}
m_{\Tc}(a\tp b)=m_{\Kc}\bigl(\Ru_1\tr a \tl \Ru_{1'}^{-1}\tp b \tl \Ru_{2'}^{-1}\tl \Ru_2\bigr).
\ee
Let $G$ be the orthogonal or symplectic groups $G$.
Formula (\ref{tw_pr}) applied  to the equation
\be
T^tB^{-1}T=B^{-1},
\quad TBT^t=B
\ee
gives
\be
\label{orhtog_re}
R_1^tK^t\bigr ( (R^t_{1'})^{-1}B^{-1} ( R_{2'})^{-1}\bigl)R_2 K= B^{-1},
\quad
K R_1 B K^tR_2^t =R_{1'}B R_{2'}^t.
\ee
The ideal generated by (\ref{orhtog_re}) lies in the kernel of  the
$\U_\hbar(\D\g)$-equivariant projection $\Kc\to \C_\hbar[G_{STS}]$.

Similarly one can express the element $\det_q(T)$ through the generators $K^i_j$.
We denote by $\det_q(K)$ the resulting form of degree $n$. The ideal $\bigl(\det_q(K)-1\bigr)$ is annihilated by
the $\U_\hbar(\D\g)$-equivariant projection $\Kc\to \C_\hbar[G_{STS}]$.
\begin{propn}
\label{qSTS}
Let $G$ be a classical complex matrix group. Then
the algebra $\C_\hbar[G_{STS}]$ is isomorphic to
the quotient of $\Kc$ by the ideal $\Jg$, where

i) $\Jg$ is generated by relations (\ref{orhtog_re}) and $\det_q(K)=1$
for the   $G=SO(N)$.

ii) $\Jg=\bigl(\det_q(K)-1\bigr)$ for $G=SL(n)$.

iii) $\Jg$ is generated by relations (\ref{orhtog_re}) for $G=O(N)$ or $G=Sp(n)$.
\end{propn}
\begin{proof}
This proposition can be proven by a straightforward modification of the proof of Proposition \ref{sharp}.
Another way is to start from Proposition \ref{sharp} and use the RE twist applied to the DS quantization,
cf. Appendix \ref{apA}.
\end{proof}
\appendix
\section{On twist of module algebras}
\label{apA}
In this subsection we study how twist of Hopf algebras affects defining relations
of their module algebras.
Let $\Ha$ be a Hopf algebra, $V$ a finite dimensional left $\Ha$-module,
and $T(V)$ the tensor algebra of $V$. Let $W$ be an $\Ha$-submodule in $T(V)$
generating an ideal $J(W)$ in $T(V)$. Denote by $\A$ the quotient algebra $T(V)/J(W)$.

Let $\F\in \Ha\tp \Ha$ be a twisting cocycle and $\tilde \Ha$ the corresponding twist of $\Ha$.
Denote by $\tilde \A$ the twist of the module algebra $\A$.
The multiplication in $\tilde\A$ is expressed through the multiplication in $\A$ by
$m_{\tilde \A}=m_{\A}\circ \F$ and similarly for $\widetilde {T(V)}$ and ${T(V)}$.

For each $n=0,1,\ldots,$ introduce an automorphism of $V^{\tp n}$ by induction:
$$\Omega_n=\id, \quad n=0,1, \quad \Omega_n=(\Delta^m\tp \Delta^k)(\F)\bigl(\Omega_m\tp \Omega_k\bigr),\quad k+m=n.$$
This definition does not depend on the partition $k+m=n$.
The elements $\Omega_n$ amounts to a linear automorphism $\Omega$ of $T(V)$.
\begin{propn}
The algebra $\tilde\A$ is isomorphic to the quotient algebra $T(V)/{J(\Omega^{-1}W)}$.
\label{propn_ap}
\end{propn}
\begin{proof}
Since the ideal $J(W)\subset T(V)$ is invariant, it is also an ideal in $\widetilde {T(V)}$.
It is easy to see that the quotient $\widetilde {T(V)}/J(W)$ is isomorphic to  $\tilde\A$.
On the other hand, the algebra $\widetilde {T(V)}$ is isomorphic to $T(V)$. The isomorphism
is given by the maps $T(V)\supset m(v_1\tp \ldots \tp v_n)\mapsto (\tilde m\circ\Omega_n)(v_1\tp \ldots \tp v_n)$,
$n\in \N$,
where $m$ and $\tilde m$ are multiplications in $T(V)$ and $\widetilde {T(V)}$.
This implies the proposition.
\end{proof}

\end{document}